\RequirePackage{fix-cm}
\documentclass[smallextended]{svjour3}       
\smartqed  
\usepackage{graphicx}
\usepackage{bbm}
\usepackage{graphicx}

\usepackage{amssymb}
\usepackage{amsfonts}
\usepackage{mathrsfs}
\usepackage{amssymb}
\usepackage{color}
\usepackage{amsmath}
\numberwithin{equation}{section}

\numberwithin{lemma}{section}

\numberwithin{theorem}{section}

\numberwithin{proposition}{section}

\numberwithin{remark}{section}

\allowdisplaybreaks

\begin{document}

\title{Weak order in  averaging principle for stochastic differential equations with jumps
}


\author{Bengong Zhang\and Hongbo Fu  \and
        Li Wan \and Jicheng Liu}


\institute{Bengong Zhang \at
             College of Mathematics and Computer Science, Wuhan Textile University \\
              \email{benyan1219@126.com}
            \and
              Hongbo Fu \at
              College of Mathematics and Computer Science, Wuhan Textile University \\
              \email{hbfu@wtu.edu.cn}           
           \and
           Li Wan \at
              College of Mathematics and Computer Science, Wuhan Textile
              University\\
              \email{wanlinju@aliyun.com}
           \and
           Jicheng Liu \at
           School of Mathematics and Statistics, Huazhong University of Science and Technology \\
           \email{jcliu@hust.edu.cn}}

\date{Received: date / Accepted: date}

\maketitle

\begin{abstract}
The present article  deals with the  averaging principle  for a
two-time-scale system of jump-diffusion stochastic differential
equation. Under suitable conditions, the weak error is expanded in
powers of timescale parameter. It is proved that the rate of weak
convergence to the averaged  dynamics is of order $1$. This reveals
the rate of weak convergence is essentially twice that of strong
convergence.

\end{abstract}

\section{Introduction}
 We consider a two-time-scale system  of jump-diffusion stochastic
differential equation in form of
\begin{eqnarray}
&&\!\!\!\!\!\!\!\!\!\!\!\!dX^\epsilon_t=a(X_t^\epsilon,
Y_t^\epsilon)dt+b(X_t^\epsilon)d {B}_t+c(X_{t-}^\epsilon)d {P}_t,
\;X_0^\epsilon=x,\label{slow-equation}\\
&&\!\!\!\!\!\!\!\!\!\!\!\!dY^\epsilon_t=\frac{1}{\epsilon}f(X_t^\epsilon,
Y_t^\epsilon)dt+\frac{1}{\sqrt{\epsilon}}g(X_t^\epsilon,
Y_t^\epsilon)d {W}_t+h(X_{t-}^\epsilon,Y_{t-}^\epsilon)d
{N}^\epsilon_t, \;Y_0^\epsilon=y,\label{fast-equation}
\end{eqnarray}
where  $X_t^\epsilon\in \mathbb{R}^n, Y_t^\epsilon\in \mathbb{R}^m$,
the drift  functions $a(x, y)\in  \mathbb{R}^n, f(x,y)\in
\mathbb{R}^m,$ the diffusion functions $b(x )\in\mathbb{R}^{n\times
d_1}, c(x)\in \mathbb{R}^n, g(x,y)\in\mathbb{R}^{m\times d_2}$ and
$h(x,y)\in\mathbb{R}^{m}$. $B_t$ and $W_t$ are the vectors of $d_1,
d_2$-dimensional independent Brownian motions on a complete
stochastic base $(\Omega, \mathcal{F},\mathcal{F}_t, \mathbb{P})$,
respectively. $P_t$ is a scalar simple Poisson process with
intensity $\lambda_1$, and $N_t^\epsilon$ is a scalar simple Poisson
process with intensity $\frac{\lambda_2}{\epsilon}$. The positive
parameter $\epsilon$ is small and  describes the ratio of time
scales between $X^\epsilon_t$ and $Y^\epsilon_t$. Systems
\eqref{slow-equation}-\eqref{fast-equation} with two time scales
occur frequently in applications  including chemical kinetics,
signal processing, complex fluids and  financial engineering.

With the  separation of time scale,  we can view the state variable
of the system as being divided into two parts: the ``slow" variable
$X^\epsilon_t$ and the ``fast" variable $Y^\epsilon_t$. It is often
the case that we are interested only in the dynamics of slow
component. Then a simplified equation, which is independent of fast
variable and possesses the essential features of the system, is
highly desirable. Such a simplified equation is often constructed by
averaging procedure as in \cite{Bogoliubov,Volosov} for
deterministic ordinary differential equations, as well as the
further development
\cite{Khas,Freidlin-Wentzell1,Freidlin-Wentzell2,Vere1,Vere2,Kifer1,Kifer2,Kifer3,Wainrib}
for stochastic differential equations with continuous Gaussian
processes. As far as averaging for stochastic dynamical systems in
infinite dimensional space is concerned, it is worthwhile to quote
the important works of \cite{Cerrai1,Cerrai2,Cerrai-Siam,wangwei}
and also the works of \cite{Fu-Liu,Fu-Liu-2,Xujie}. For related
works on averaging for multivalued stochastic differential equations
we refer the reader to \cite{Guo,Xujie2}.

In order to derive the averaged dynamics of the system
\eqref{slow-equation}-\eqref{fast-equation}, we introduce the fast
motion equation with  a frozen slow component $x\in \mathbb{R}^n$ in
form of
\begin{eqnarray}\label{frozen}
dY_t^x= f(x, Y_t^x)dt+g(x, Y_t^x)d {W}_t+h(x,Y^x_{t-})d {N}_t,
\;Y_0^x=y,
\end{eqnarray}
whose solution is denoted by $Y_t^\epsilon(y)$. Under suitable
conditions on $f,g$ and $h$, $Y_t^\epsilon(y)$ induces a unique
invariant  measure $\mu^x(dy)$ on $\mathbb{R}^m$, which   is ergodic
and ensures the averaged equation:
\begin{eqnarray*}
d\bar{X}_t=\bar{a}(\bar{X}_t)dt+b(\bar{X}_t)d{B}_t+c(\bar{X}_{t-})d{P}_t,
\;\bar{X}_0 =x,
\end{eqnarray*}
where the averaging nonlinearity is defined by setting
\begin{eqnarray*}
\bar{a}(x )&=&\int_{\mathbb{R}^m}a(x,y)\mu^x(dy)\\
&=&\lim_{t\rightarrow +\infty} \mathbb{E} a(x, Y_t^x(y)).
 \end{eqnarray*}

In \cite{Givon}, it was shown that under the above conditions the
slow motion $X^\epsilon_t$ converges strongly to the solution
$\bar{X}_t$ of the above averaged  equation with jumps. The order of
convergence $\frac{1}{2}$ in strong sense was provided in
\cite{LiuDi1}. { To our best knowledge, there is no existing
literature to address the weak order in averaging principle for jump
diffusion stochastic differential systems. In fact, it is fair to
say that the weak convergence in stochastic averaging theory of
systems driven by jump noise is not  fully developed yet, although
some strong approximation results on the rate of strong convergence
were obtained \cite{Bao,Xu,Xu3}.

Therefore, we aim to study this problem in this paper. Here we are
interested in  the rate of weak convergence of the averaging
dynamics to the true solution of slow motion $X^\epsilon_t$. In
other word,  we will determine the order, with respect to timescale
parameter $\epsilon$, of weak deviation between original solution to
slow equation and the solution  of the corresponding averaged
equation.} The main technique we adapted is to find an expansion
with respect to $\epsilon$ of the solutions of the Kolmogorov
equations associated with the jump diffusion system. The solvability
of the Poisson equation associated with the generator of frozen
equation provides an expression for the coefficients of the
expansion. As a result, the boundedness for the coefficients of
expansion can be proved by smoothing effect of the corresponding
transition semigroup in the space of bounded and uniformly
continuous functions, where some regular conditions is needed on
drift and diffusion term.

 { Our result shows that the weak
convergence rate to be $1$  even when there are jump components in
the system. It is the main contribution of this work.} We would like
to stress that asymptotic  method was first applied by Br\'{e}hier
\cite{Brehier} to an averaging result for stochastic
reaction-diffusion equations in the case of Gaussian noise of
additive type was included only in the fast motion. However, the
extension of this argument  is not straightforward. The method used
in the proof of weak order in \cite{Brehier} is strictly related to
the differentiability in time of averaged process. Therefore, once
the noise is introduced in the slow equation, difficulties will
arise and the procedure becomes more complicated. Our result in this
paper bridges such a gap, in which the slow and the fast motions are
both perturbed  by  noise with jumps.

The rest of the paper is  structured as follows.  Section 2 is
devoted to  notations, assumptions and summarize preliminary
results. The ergodicity of fast process and the averaged dynamics of
system with jumps is introduced in Section 3. Then the main result
of this article, which is derived via the asymptotic expansions and
uniform error estimates, is presented in Section 4. Finally, we give
the appendix in  section 5.

It should be pointed out that the letter $C$ below with or without
subscripts will denote  generic positive constants independent of
$\epsilon$ in the whole paper.

\section{Assumptions and preliminary results}
For any integer $d$, the scalar product and norm on $d-$dimensional
Euclidean space $\mathbb{R}^d$ are denoted by
$\Big(\cdot,\cdot\Big)_{\mathbb{R}^d}$ and
$\|\cdot\|_{\mathbb{R}^d}$, respectively. For any  integer $k$, we
denote by $C_b^k(\mathbb{R}^d,\mathbb{R})$ the space of all
$k-$times differentiable functions on $\mathbb{R}^d$, which have
bounded uniformly continuous derivatives up to the $k$-th order.

In what follows, we shall assume that the drift and diffusion
coefficients arising in the system  fulfill  the following
conditions.

 (\textbf{A1})
The mappings $a(x,y), b(x), c(x), f(x,y), g(x,y)$ and $h(x,y)$ are
of class $C^2$ and have bounded first and second derivatives.
Moreover, we assume that $a(x,y), b(x)$ and $c(x)$ are bounded.

(\textbf{A2}) There exists a constant $\alpha>0$ such that for any
$x\in \mathbb{R}^n, y\in\mathbb{R}^m$ it holds
\begin{eqnarray*}
y^Tg(x,y)g^T(x,y)y\geq \alpha\|y\|_{\mathbb{R}^m}.
\end{eqnarray*}

(\textbf{A3}) There exists a constant $\beta>0$ such that for any
$y_1, y_2\in\mathbb{R}^m$ and $x\in\mathbb{R}^n$ it holds
\begin{eqnarray*}
&&\Big(y_1-y_2,
f(x,y_1)-f(x,y_2)+\lambda_2(h(x,y_1)-h(x,y_2))\Big)_{\mathbb{R}^m}\\
&&+\|g(x,y_1)-g(x,y_2)\|^2_{\mathbb{R}^m}+\lambda_2|h(x,y_1)-h(x,y_2)|^2\\
&&\leq -\beta \|y_1-y_2\|^2_{\mathbb{R}^m}.
\end{eqnarray*}

\begin{remark}
Notice that from (\textbf{A1}) it immediately follows that the
following directional derivatives exist and are controlled:
\begin{eqnarray*}
&&\|D_xa(x,y)\cdot k_1\|_{\mathbb{R}^n}\leq
L\|k_1\|_{\mathbb{R}^n},\\
&&\|D_ya(x,y)\cdot l_1\|_{\mathbb{R}^n}\leq
L\|l_1\|_{\mathbb{R}^m},\\
&&\|D_{xx}^2a(x,y)\cdot(k_1,k_2)\|_{\mathbb{R}^n}\leq
L\|k_1\|_{\mathbb{R}^n}\|k_2\|_{\mathbb{R}^n},\\
&&\|D_{yy}^2a(x,y)\cdot(l_1,l_2)\|_{\mathbb{R}^n}\leq
L\|l_1\|_{\mathbb{R}^m}\|l_2\|_{\mathbb{R}^m},
\end{eqnarray*}
where $L$ is a constant independent of $x,y,k_1, k_2, l_1$ and
$l_2$. For differentiability of mappings $b, c,f, g$ and $h$ we
possess the analogous  results. For examples,  we have
\begin{eqnarray*}
&&\|D^2_{xx}b(x)\cdot (k_1,k_2)\|_{\mathbb{R}^n}\leq L
\|k_1\|_{\mathbb{R}^n} \|k_2\|_{\mathbb{R}^n},\; k_1,k_2\in \mathbb{R}^n,\\
&&\|D^2_{yy}f(x,y)\cdot (l_1,l_2)\|_{\mathbb{R}^m}\leq L
\|l_1\|_{\mathbb{R}^m} \|l_2\|_{\mathbb{R}^m}\|,\;l_1,l_2\in
\mathbb{R}^m.
\end{eqnarray*}

As far as the assumption (\textbf{A2}) is concerned, it is a sort of
non-degeneracy condition and it is assumed in order to have the
regularizing effect of the Markov transition semigroup associated
with   the fast dynamics. Assumption (\textbf{A3}) is the
dissipative condition   which determines how the fast equation
converges to its equilibrium state.
\end{remark}

As assumption (\textbf{A1}) holds,  for any $\epsilon>0$ and any
initial conditions $x\in \mathbb{R}^n$ and $y\in \mathbb{R}^m$,
system \eqref{slow-equation}-\eqref{fast-equation} admits a unique
solution, which, in order to emphasize the dependence on the initial
data, is denoted by $(X_t^\epsilon(x,y), Y_t^\epsilon(x,y))$.
Moreover the following lemma holds (for a proof see e.g.
\cite{LiuDi1}).
\begin{lemma}
Under the assumptions (\textbf{A1}), (\textbf{A2}) and
(\textbf{A3}),
 for any $x\in \mathbb{R}^n$, $y\in \mathbb{R}^m$ and $\epsilon>0$
 we have
\begin{eqnarray} \mathbb{E}\|X_t^\epsilon(x,y)\|^2_{\mathbb{R}^n}\leq
C_T(1+\|x\|^2_{\mathbb{R}^n}+\|y\|^2_{\mathbb{R}^m}), \;t\in [0,
T]\label{X-bound}
\end{eqnarray}
and
\begin{eqnarray}
\mathbb{E}\|Y_t^\epsilon(x,y)\|^2_{\mathbb{R}^n}\leq
C_T(1+\|x\|^2_{\mathbb{R}^m}+\|y\|^2_{\mathbb{R}^m}), \;t\in [0,
T].\label{Y-bound}
\end{eqnarray}
\end{lemma}

\section{Frozen equation and averaged equation }
Fixing $\epsilon=1$, we consider the fast equation with frozen slow
component $x\in\mathbb{R}^n$,
\begin{eqnarray}\label{frozen}
\begin{cases}
dY_t^x(y)= f(x, Y_t^x(y))dt+g(x, Y_t^x(y))d {W}_t+h(x,Y^x_{t-}(y))d
{N}_t, \\
Y_0^x=y.
\end{cases}
\end{eqnarray}
Under assumptions (\textbf{A1})-(\textbf{A3}), such a problem has a
unique solution, which satisfies \cite{LiuDi1}:
\begin{eqnarray}\label{frozen-bound}
\mathbb{E}\|Y_t^x(y)\|^2_{\mathbb{R}^m}\leq
C(1+\|x\|^2_{\mathbb{R}^n}+e^{-\beta
t}\|y\|^2_{\mathbb{R}^m}),\;t\geq 0.
\end{eqnarray}
Let $Y_t^{x}( y')$ be the solution of problem \eqref{frozen} with
initial value $Y_0^x=y'$, the It\^{o} formula implies that for any
$t\geq0$,
\begin{eqnarray}
\mathbb{E}\|Y_t^{x}(y)-Y_t^{x}
(y')\|_{\mathbb{R}^m}^2\leq\|y-y'\|_{\mathbb{R}^m}^2e^{-\beta t}.
\label{difference-initial}
\end{eqnarray}
Moreover, as discussion in \cite {LiuDi1}  and \cite{Givon},
equation \eqref{frozen} admits a unique ergodic invariant measure
$\mu^x$ satisfying
\begin{equation}\label{mu-Momenent-bound}
\int_{\mathbb{R}^m}\|y\|_{\mathbb{R}^m}^2\mu^x(dy)\leq
C(1+\|x\|^2_{\mathbb{R}^n}).
\end{equation}
Then, by averaging the coefficient $a$ with respect to the invariant
measure $\mu^x$, we can define an $\mathbb{R}^n$-valued mapping
\begin{equation*}
\bar{a}(x):=\int_{\mathbb{R}^m}a(x,y)\mu^x(dy), x\in \mathbb{R}^n.
\label{aver-a}
\end{equation*}
Due to assumption (\textbf{A1}), it is easily to check that
$\bar{a}(x)$ is 2-times differentiable with bounded derivatives, and
hence it is Lipschitz-continuous such that
\begin{eqnarray*}
\|\bar{a}(x_1)-\bar{a}(x_2)\|_{\mathbb{R}^n}\leq
C\|x_1-x_2\|_{\mathbb{R}^n},\;x_1, x_2\in \mathbb{R}^n.
\label{bar-a-lip}
\end{eqnarray*}
According to   invariant property of $\mu^x$,
\eqref{mu-Momenent-bound} and assumption (\textbf{A1}), we have
\begin{eqnarray}
\nonumber\left\|\mathbb{E}a(x, Y_t^{x}(
y))-\bar{a}(x)\right\|^2_{\mathbb{R}^n}&=&
\|\int_{\mathbb{R}^m}\mathbb{E}\big(a(x, Y_t^{x}
(y))-a(x, Y_t^{x}(z))\big)\mu^x(dz) \|^2_{\mathbb{R}^n}\\
\nonumber&\leq&\int_{\mathbb{R}^m}\mathbb{E}\left\|Y_t^{x}
(y)-Y_t^{x}(z)\right\|_{\mathbb{R}^m}^2\mu^x(dz)\\
\nonumber&\leq&e^{-\beta t}\int_{\mathbb{R}^m}\|y-z\|^2_{\mathbb{R}^m}\mu^x(dz)\\
&\leq&Ce^{-\beta
t}\big(1+\|x\|^2_{\mathbb{R}^n}+\|y\|^2_{\mathbb{R}^m}\big).\label{Averaging-Expectation}
\end{eqnarray}
Now we can  introduce the effective dynamical system
\begin{eqnarray}\label{averaging-equation}
\begin{cases}
d\bar{X}_t(x)=\bar{a}(\bar{X}_t(x))dt+b(\bar{X}_t(x))d{B}_t+c(\bar{X}_{t-}(x))d
{P}_t,
\\
\bar{X}_0 =x.
\end{cases}
\end{eqnarray}
As the coefficients  $\bar{a}, b$ and $c$ are Lipschitz-continuous,
this equation  admits a unique solution such that
\begin{eqnarray}\label{bar-x-bound}
\mathbb{E}\|\bar{X}_t(x)\|^2_{\mathbb{R}^n}\leq
C_T(1+\|x\|^2_{\mathbb{R}^n}), \;t\in [0, T].
\end{eqnarray}

 With the above assumptions and notations we have the following   result, which
is a direct consequence of Lemma \ref{u-0}, Lemma
\ref{u-1-abso-lemma} and Lemma \ref{r-bound-lemma}.

\begin{theorem}
Assume that $x\in \mathbb{R}^n$ and $y\in\mathbb{R}^m$, Then, under
assumptions (\textbf{A1}), (\textbf{A2}) and (\textbf{A3}), for  any
$T>0$ and $\phi\in C_b^3(\mathbb{R}^n,\mathbb{R})$, there exists a
 constant $C_{T,\phi,x,y}$ such that
 \begin{eqnarray*}
 \left|\mathbb{E}\phi(X^\epsilon_T(x,y))-\mathbb{E}\phi(\bar{X}_T(x))\right|\leq
 C_{T,\phi,x,y}\epsilon.
 \end{eqnarray*}
 As a consequence, it can be claimed that the weak order in
 averaging principle for jump-diffusion stochastic  systems is 1.
\end{theorem}

\section{Asymptotic expansion}
Let $\phi\in C_b^3(\mathbb{R}^n, \mathbb{R})$ and define a function
$u^\epsilon(t, x,y):[0, T]\times
\mathbb{R}^n\times\mathbb{R}^m\rightarrow \mathbb{R}$ by
\begin{equation*}
u^\epsilon(t, x,y)=\mathbb{E}\phi(X_t^\epsilon(x,y)).
\end{equation*}
We are now ready to seek an expansion  formula for $u^\epsilon(t,
x,y)$ with respect to $\epsilon$ with the form
\begin{eqnarray}\label{asymp-expan}
u^\epsilon(t,x,y)=u_0(t,x,y)+\epsilon u_1(t,x,y)+r^\epsilon(t,x,y),
\end{eqnarray}
where  $u_0$ and $u_1$ are smooth functions which will be
constructed below, and $r^\epsilon$ is the remainder term. To this
end, let us recall the Kolmogorov operator corresponding to the slow
motion equation, with a frozen fast component $y\in\mathbb{R}^m$,
which is a  second order   operator  taking form
\begin{eqnarray}
\mathcal{L}_1\Phi(x)&=&\Big(a(x,y),D_x
\Phi(x)\Big)_{\mathbb{R}^n}+\frac{1}{2}Tr\left[D_{xx}^2 \Phi(x)
\cdot
b(x ) b^T(x )\right]\nonumber\\
&&+\lambda_1(\Phi(x+c(x ))-\Phi(x)), \;\Phi\in C_b^2(\mathbb{R}^n,
\mathbb{R}).\nonumber
\end{eqnarray}
For any frozen slow component $x\in\mathbb{R}^m$, the Kolmogorov
operator for equation \eqref{frozen} is given by
\begin{eqnarray}
\mathcal{L}_2\Psi(y)&=&\Big(f(x,y),D_y
\Psi(y)\Big)_{\mathbb{R}^m}+\frac{1}{2}Tr\left[  D_{yy}^2\Psi(y)
\cdot
g(x,y) g^T(x,y)\right]\nonumber\\
&&+\lambda_2(\Psi(y+h(x,y))-\Psi(y)), \;\Psi\in C_b^2(\mathbb{R}^m,
\mathbb{R}).\nonumber
\end{eqnarray}
We  set
\begin{equation*}
\mathcal{L}^\epsilon:=\mathcal{L}_1+\frac{1}{\epsilon}\mathcal{L}_2.
\end{equation*}
It is known $u^\epsilon(t,x,y)$ solves the  equation
\begin{eqnarray}\label{Kolm}
\begin{cases}
\frac{\partial}{\partial t}u^\epsilon(t, x, y)=\mathcal {L}^\epsilon u^\epsilon(t, x, y),\\
u^\epsilon(0, x,y)=\phi(x),
\end{cases}
\end{eqnarray}
Also recall the Kolmogorov operator associated with the averaged
equation \eqref{averaging-equation} is defined as
\begin{eqnarray}
\bar{\mathcal{L}} \Phi(x)&=&\Big(\bar{a}(x ),D_x
\Phi(x)\Big)_{\mathbb{R}^n}+\frac{1}{2}Tr\left[  D_{xx}^2\Phi(x)
\cdot
b(x ) b^T(x )\right]\nonumber\\
&&+\lambda_1(\Phi(x+c(x))-\Phi(x)), \;\Phi\in C_b^2(\mathbb{R}^n,
\mathbb{R}).\nonumber
\end{eqnarray}
If we set $$\bar{u}(t, x)=\mathbb{E}\phi(\bar{X}_t(x)),$$ we have
\begin{eqnarray}\label{Kolm-Aver}
\begin{cases}
\frac{\partial}{\partial t}\bar{u}(t, x)=\bar{\mathcal {L}} \bar{u}(t, x),\\
\bar{u}(0, x)=\phi(x).
\end{cases}
\end{eqnarray}

\subsection{\textbf{The leading term}}
Let us begin with constructing the leading term.
 By substituting expansion \eqref{asymp-expan} into \eqref{Kolm}, we
see that
\begin{eqnarray*}
\frac{\partial u_0}{\partial t}+\epsilon\frac{\partial u_1}{\partial
t}+\frac{\partial r^\epsilon}{\partial t}&=& \mathcal
{L}_1u_0+\epsilon \mathcal {L}_1u_1+\mathcal
{L}_1r^\epsilon\\
&+&\frac{1}{\epsilon}\mathcal {L}_2u_0+\mathcal {\mathcal
{L}}_2u_1+\frac{1}{\epsilon}\mathcal {L}_2r^\epsilon.
\end{eqnarray*}
By  equating powers of      $\epsilon$, we obtain the following
system of equations:
\begin{eqnarray}
&&\mathcal {L}_2u_0=0, \label{u-o-equ-1} \\
&&\frac{\partial u_0}{\partial t}=\mathcal {L}_1u_0+\mathcal
{L}_2u_1.\label{u-0-equ-2}
\end{eqnarray}
According to \eqref{u-o-equ-1}, we can conclude $u_0$ does not
depend on $y$, that is $$u_0(t,x, y)=u_0(t,x).$$ We also impose the
initial condition $u_0(0,x)=\phi(x).$ Note that  $\mathcal {L}_2$ is
the generator of a Markov process defined by equation
\eqref{frozen}, which admits  a unique invariant measure $\mu^x$, we
have
\begin{eqnarray}
\int_{\mathbb{R}^m}\mathcal
{L}_2u_1(t,x,y)\mu^x(dy)=0.\label{u-1-zero}
\end{eqnarray}
Thanks to   \eqref{u-0-equ-2}, this yields
\begin{eqnarray*}
\frac{\partial u_0}{\partial t}(t,x)&=&\int_{\mathbb{R}^m}\frac{\partial u_0}{\partial t}(t,x)\mu^x(dy)\\
&=&\int_{\mathbb{R}^m}\mathcal {L}_1u_0(t,x)\mu^x(dy)\\
&=&\int_{\mathbb{R}^m}\Big(a(x,y), D_x
u_0(t,x)\Big)_{\mathbb{R}^n}\mu^x(dy)\\
&&+\frac{1}{2}Tr\left[D_{xx}^2u_0(t,x)\cdot
b(x ) b^T(x )\right]\nonumber\\
&&+\lambda_1(u_0(x+c(x ))-u_0(x))\\ &=&\bar{\mathcal {L}}u_0(t,x),
\end{eqnarray*}
so that $u_0$ and $\bar{u}$ are described by the same evolutionary
equation. By  uniqueness argument, we easily have the following
lemma:
\begin{lemma}\label{u-0}
Under assumptions (\textbf{A1}), (\textbf{A2}) and (\textbf{A3}),
for any $x\in \mathbb{R}^n$, $y\in \mathbb{R}^m$ and $T>0$, we have
$u_0(T,x,y)=\bar{u}(T,x)$.
\end{lemma}

\subsection{\textbf{Construction of} $ {u_1}$}
According to Lemma \ref{u-0}, \eqref{Kolm-Aver} and
\eqref{u-0-equ-2}, we get
\begin{eqnarray*}
\bar{\mathcal {L}}\bar{u}=\mathcal{L}_1\bar{u}+\mathcal{L}_2u_1,
\end{eqnarray*}
which means that
\begin{eqnarray}
\mathcal
{L}_2u_1(t,x,y)&=&\Big(\bar{a}(x)-a(x,y),D_x\bar{u}(t,x)\Big)_{\mathbb{R}^n}\nonumber\\
&:=&-\rho(t,x,y),\label{equation-u1}
\end{eqnarray}
where $\rho$ is of class $C^2$ with respect to $y$, with uniformly
bounded derivatives. Moreover, for any $t\geq 0$ and $x\in
{\mathbb{R}^n} $, the equality \eqref{u-1-zero} guarantees that
\begin{eqnarray*}
\int_{\mathbb{R}^m}\rho(t,x,y)\mu^x(dy)=0.
\end{eqnarray*}
For any $y\in \mathbb{R}^m$ and $s>0$ we have
\begin{eqnarray}
\frac{\partial}{\partial s}\mathcal
{P}_s\rho(t,x,y)&=&\Big(f(x,y),D_y [\mathcal
{P}_s\rho(t,x,y)]\Big)_{\mathbb{R}^m}\nonumber\\
&&+\frac{1}{2}Tr\left[D_{yy}^2[\mathcal {P}_s\rho(t,x,y)]\cdot
g(x,y) g^T(x,y)\right]\nonumber\\
&&+\lambda_2\left(\mathcal {P}_s[\rho(t,x,y+h(x,y))]-\mathcal
{P}_s[\rho(t,x,y)]\right),\label{P-s-rho}
\end{eqnarray}
here $$\mathcal {P}_s[\rho(t,x,y)]:=\mathbb{E}\rho(t, x,Y^x_s(y)).$$
Recalling that   $\mu^x$ is the unique invariant measure
corresponding to Markov process $Y^x_t(y)$ defined by  equation
\eqref{frozen}, from Lemma \ref{bar-u-x} we infer that
\begin{eqnarray*}
&&\left|\mathbb{E}\rho(t,
x,Y^x_s(y))-\int_{\mathbb{R}^m}\rho(t,x,z)\mu^x(dz)\right|\nonumber\\
&&=\left|\int_{\mathbb{R}^m}\mathbb{E}[\rho(t, x,Y^x_s(y))- \rho(t,
x,Y^x_s(z))]\mu^x(dz)\right|\nonumber\\
&&\leq\int_{\mathbb{R}^m}\left|\mathbb{E}\Big(a(x,Y^x_s(z))-
a(x,Y^x_s(y)), D_x\bar{u}(t,x)\Big)_{\mathbb{R}^n}\right|\mu^x(dz)\nonumber\\
&&\leq C\int_{\mathbb{R}^m} \mathbb{E} \|Y_s^x(z)-
Y_s^x(y)\|_{\mathbb{R}^n}
\mu^x(dz).\nonumber\\
\end{eqnarray*}
Now it follows from   \eqref{difference-initial} and
\eqref{mu-Momenent-bound} that
\begin{eqnarray*}
&&\left|\mathbb{E}\rho(t,
x,Y_s^x(y))-\int_{\mathbb{R}^m}\rho(t,x,z)\mu^x(dz)\right|\nonumber\\
&&\leq
C(1+\|x\|_{\mathbb{R}^n}+\|y\|_{\mathbb{R}^m})e^{-\frac{\beta}{2}s},
\end{eqnarray*}
which implies
\begin{equation}
\lim\limits_{s\rightarrow{+\infty}}\mathbb{E}\rho(t,
x,Y_s^x(y))=\int_{\mathbb{R}^m}\rho(t,x,z)\mu^x(dz)=0. \nonumber
\end{equation}
With the aid of the above limit, we can deduce from \eqref{P-s-rho}
that
\begin{eqnarray*}
&&\Big(f(x,y),D_y \int_0^{+\infty}[\mathcal
{P}_s\rho(t,x,y)]ds\Big)_{\mathbb{R}^m}\\
&&+\frac{1}{2}Tr\left[D_{yy}^2\int_0^{+\infty}[\mathcal
{P}_s\rho(t,x,y)]\cdot
g(x,y) g^T(x,y)ds\right]\nonumber\\
&&+\lambda_2\left(\int_0^{+\infty}\mathcal
{P}_s[\rho(t,x,y+h(x,y))]ds-\int_0^{+\infty}\mathcal
{P}_s[\rho(t,x,y)]ds\right)\nonumber\\
&&=\int_0^{+\infty}\frac{\partial}{\partial s}\mathcal{P}_s[\rho(t,x,y)]ds\\
&&=\lim\limits_{s\rightarrow{+\infty}}\mathbb{E}\rho(t,
x,Y_s^x(y))-\rho(t,x,y)\nonumber\\
&&=\int_{\mathbb{R}^m}\rho(t,x,z)\mu^x(dz)-\rho(t,x,y)\nonumber\\
&&=-\rho(t,x,y),\nonumber
\end{eqnarray*}
which implies $$\mathcal{L}_2(\int_0^{{+\infty}}\mathcal
{P}_s\rho(t,x,y) ds)=-\rho(t,x,y).$$ Therefore,
\begin{eqnarray}\label{u-1}
u_1(t,x,y):=\int_0^{+\infty} \mathbb{E}\rho(t,x,Y^x_s(y))ds
\end{eqnarray}
is the solution to  equation \eqref{equation-u1}.
\begin{lemma}\label{u-1-abso-lemma}
Under assumptions (\textbf{A1}), (\textbf{A2}) and (\textbf{A3}),
 for any $x\in \mathbb{R}^n$, $y\in \mathbb{R}^m$ and $T>0$, we have
\begin{eqnarray}
|u_1(t,x,y)|\leq C_T(1
+\|x\|_{\mathbb{R}^n}+\|y\|_{\mathbb{R}^m}),\;t\in[0, T].
\label{u-1-abso}
\end{eqnarray}
\begin{proof}
By \eqref{u-1}, we have
\begin{equation*}
u_1(t,x,y)=\int_0^{{+\infty}}\mathbb{E}\Big(\bar{a}(x)- \
a(x,Y_s^x(y)), D_x\bar{u}(t,x)\Big)_{\mathbb{R}^n }ds,
\end{equation*}
so that
\begin{eqnarray*}
|u_1(t,x,y)|&\leq&\int_0^{{+\infty}}\|\bar{a}(x)-
\mathbb{E}[a(x,Y_s^x(y))]\|_{\mathbb{R}^n }
\cdot\|D_x\bar{u}(t,x)\|_{\mathbb{R}^n } ds.
\end{eqnarray*}
Therefore, from Lemma \ref{bar-u-x} and
\eqref{Averaging-Expectation}, we get
\begin{eqnarray*}
 |u_1(t,x,y)|&\leq& C_T(1+\|x\|_{\mathbb{R}^n }  +\|y\|_{\mathbb{R}^m })\int_0^{{+\infty}}e^{-\frac{\beta}{2} s}ds\\
&\leq&C_T(1+\|x\|_{\mathbb{R}^n }  +\|y\|_{\mathbb{R}^m }).
\end{eqnarray*}
\end{proof}
\end{lemma}

\subsection{\textbf{Determination of remainder} $ {r^\epsilon}$}
We now turn to  the construction for   remainder term $r^\epsilon$.
It is known that
\begin{eqnarray*}
(\partial_t-\mathcal{L}^\epsilon)u^\epsilon=0,
\end{eqnarray*}
which, together with \eqref{u-o-equ-1} and \eqref{u-0-equ-2},
implies
\begin{eqnarray}
(\partial_t-\mathcal{L}^\epsilon)r^\epsilon&=&-(\partial_t-\mathcal{L}^\epsilon)
u_0-\epsilon(\partial_t-\mathcal{L}^\epsilon)u_1\nonumber\\
&=&-(\partial_t-\frac{1}{\epsilon}\mathcal{L}_2-\mathcal{L}_1)u_0-\epsilon(\partial_t-\frac{1}{\epsilon}\mathcal{L}_2-\mathcal{L}_1)u_1\nonumber\\
&=&\epsilon(\mathcal{L}_1u_1-\partial_tu_1).  \label{r-equation}
\end{eqnarray}
In order to estimate the remainder term $r^\epsilon$  we need the
following   two lemmas.

\begin{lemma}\label{4.3}
Under assumptions (\textbf{A1}), (\textbf{A2}) and (\textbf{A3}),
for any $x\in \mathbb{R}^n$, $y\in \mathbb{R}^m$ and $T>0$, we have
\begin{eqnarray*}
\left|\frac{\partial u_1}{\partial t}(t,x,y)\right|\leq
C_T(1+\|x\|_{\mathbb{R}^n}+\|y\|_{\mathbb{R}^m}).
\end{eqnarray*}
\begin{proof}
In view of \eqref{u-1}, we get
\begin{equation*}
\frac{\partial u_1}{\partial
t}(t,x,y)=\int_0^{{+\infty}}\mathbb{E}\left(\bar{a}(x)-
a(x,Y_s^x(y)), \frac{\partial}{\partial
t}D_x\bar{u}(t,x)\right)_{\mathbb{R}^n}ds.
\end{equation*}
By Lemma \ref{deri-xt-bar-u} introduced in Section \ref{appendix},
we have
\begin{eqnarray*}
\left|\frac{\partial u_1}{\partial
t}(t,x,y)\right|&\leq&\int_0^{{+\infty}}\mathbb{E}\left(\|\bar{a}(x)-
a(x,Y_s^x(y))\|_{\mathbb{R}^n}\cdot
\|\frac{\partial}{\partial t}D_x\bar{u}(t,x)\|_{\mathbb{R}^n}\right)ds\\
&\leq&C_T\int_0^{{+\infty}}\mathbb{E}\|\bar{a}(x)-
a(x,Y_s^x(y))\|_{\mathbb{R}^n}ds,
\end{eqnarray*}
so that from  \eqref{Averaging-Expectation} we have
\begin{eqnarray*}
\left|\frac{\partial u_1}{\partial t}(t,x,y)\right|\leq
 C_T(1+\|x\|_{\mathbb{R}^n}+\|y\|_{\mathbb{R}^m}).
\end{eqnarray*}
\end{proof}
\end{lemma}
\begin{lemma}\label{L2u1}
Under assumptions (\textbf{A1}), (\textbf{A2}) and (\textbf{A3}),
for any $x\in \mathbb{R}^n$, $y\in \mathbb{R}^m$ and $T>0$, we have
\begin{eqnarray*}
\left|\mathcal {L}_1u_1(t,x,y)\right|\leq
C_T(1+\|x\|_{\mathbb{R}^n}+\|y\|_{\mathbb{R}^m}), \;t\in [0, T].
\end{eqnarray*}
\begin{proof}
Recalling that  $u_1(t,x,y)$ is the solution of equation
\eqref{equation-u1} and equality \eqref{u-1} holds, we have
\begin{eqnarray}
\mathcal{L}_1u_1(t,x,y)&=&\Big(a(x,y),D_x
u_1(t,x,y)\Big)_{\mathbb{R}^n}\nonumber\\
&+&\frac{1}{2}Tr\left[D^2_{xx}u_1(t,x,y)\cdot
b(x ) b^T(x )\right]\nonumber\\
&+&\lambda_1[u_1(t,x+c(x ),y)-u_1(t,x,y)],\label{L1u1}
\end{eqnarray}
and then, in order to prove the boundedness of $\mathcal{L}_1u_1$,
we have to estimate the three terms arising in the right hand side
of above  equality.\\
\textbf{Step 1:}  Estimate of $\Big(a(x,y),D_x
u_1(t,x,y)\Big)_{\mathbb{R}^n}$.\\
For any $k\in \mathbb{R}^n$, we have
\begin{eqnarray*}
D_xu_1(t,x,y)\cdot k&=&\int_0^{{+\infty}}\Big(D_x\left(\bar{ a}(x)-
\mathbb{E}a(x,Y_s^x(y))\right)\cdot k,D_x\bar{u}(t,x)\Big)_{\mathbb{R}^n }ds\\
&&+\int_0^{{+\infty}}\Big(\bar{a}(x)-\mathbb{E}
a(x,Y_s^x(y)),D^2_{xx}\bar{u}(t,x)\cdot
k\Big)_{\mathbb{R}^n }ds\\
&=:&I_1(t,x,y,k)+I_2(t,x,y,k).
\end{eqnarray*}
By Lemma \ref{bar-u-x} and \ref{mix-derivative}, we infer that
\begin{eqnarray}
|I_1(t,x,y,k)| &\leq&
\|D_x\bar{u}(t,x)\|_{\mathbb{R}^n}\int_0^{+\infty}\|D_x(\bar{ a}(x)-
\mathbb{E}a(x,Y_s^x(y)))\cdot k\|_{\mathbb{R}^n}ds\nonumber\\
&\leq&C_T\|k\|_{\mathbb{R}^n}(1+\|x\|_{\mathbb{R}^n}+\|y\|_{\mathbb{R}^m})\int_0^{+\infty}e^{-\frac{\beta}{2}s}ds\nonumber\\
&\leq&C_T\|k\|_{\mathbb{R}^n}(1+\|x\|_{\mathbb{R}^n}+\|y\|_{\mathbb{R}^m}).\label{I-1}
\end{eqnarray}
By Lemma \ref{bar-u-xx} and inequality
\eqref{Averaging-Expectation}, we obtain
\begin{eqnarray*}
|I_2(t,x,y,k)| &\leq&
C_T\|k\|_{\mathbb{R}^n}\int_0^{+\infty}\|\bar{a}(x)-\mathbb{E}
a(x,Y_s^x(y))\|_{\mathbb{R}^n}ds\nonumber\\
&\leq&C_T\|k\|_{\mathbb{R}^n}(1+\|x\|_{\mathbb{R}^n}+\|y\|_{\mathbb{R}^m})
\int_0^{+\infty}e^{-\frac{\beta}{2}s}ds\nonumber\\
&\leq&C_T\|k\|_{\mathbb{R}^n}(1+\|x\|_{\mathbb{R}^n}+\|y\|_{\mathbb{R}^m}).
\end{eqnarray*}
This, together with \eqref{I-1}, implies
\begin{eqnarray*}
\|D_xu_1(t,x,y)\cdot k\|\leq
C_T\|k\|_{\mathbb{R}^n}(1+\|x\|_{\mathbb{R}^n}+\|y\|_{\mathbb{R}^m}),
\end{eqnarray*}
 and then, as $a(x,y)$ is bounded, it follows
\begin{eqnarray*}
\left|\Big(a(x,y),D_x u_1(t,x,y)\Big)_{\mathbb{R}^n}\right|&\leq&
C_T\|a(x,y)\|_{\mathbb{R}^n}(1+\|x\|_{\mathbb{R}^n}+\|y\|_{\mathbb{R}^m})\nonumber\\
&\leq& C_T
(1+\|x\|_{\mathbb{R}^n}+\|y\|_{\mathbb{R}^m}).\label{step1}
\end{eqnarray*}
\textbf{Step 2:} Estimate of $Tr\left[D^2_{xx}u_1(t,x,y)\cdot b(x )
b^T(x )\right]$.\\
Since  $u_1(t,x,y)$ is given by the  representation formula
\eqref{u-1},     for any $k_1,k_2\in \mathbb{R}^n$ we have
\begin{eqnarray*}
&&D^2_{xx}u_1(t,x,y)\cdot (k_1, k_2)\\
&&=\int_0^{{+\infty}}\mathbb{E}\Big(D^2_{xx}(\bar{a}(x)-
a(x,Y_s^x(y)))\cdot (k_1, k_2),D_x\bar{u}(t,x)\Big)_{\mathbb{R}^n }ds\\
&&+\int_0^{{+\infty}}\mathbb{E}\Big(D_x(\bar{ {a}}(x)-
a(x,Y_s^x(y)))\cdot k_1,D^2_{xx}\bar{u}(t,x)\cdot
k_2\Big)_{\mathbb{R}^n }ds\\
&&+\int_0^{{+\infty}}\mathbb{E}\Big(D_x(\bar{a}(x)-
{a}(x,Y_s^x(y)))\cdot k_2,D^2_{xx}\bar{u}(t,x)\cdot
k_1\Big)_{\mathbb{R}^n }ds\\
&&+\int_0^{{+\infty}}\mathbb{E}\Big(\bar{a}(x)- a(x,Y_s^x(y)),
D^3_{xxx}\bar{u}(t,x)\cdot (k_1,k_2)\Big)_{\mathbb{R}^n
}ds\\
&&:=\sum\limits_{i=1}^4J_i(t,x,y,k_1,k_2).
\end{eqnarray*}
Thanks to  Lemma \ref{bar-u-x} and Lemma \ref{mix-derivative-2} we
get
\begin{eqnarray}
&&\!\!\!\!\!\!\!\!\!\!|J_1(t,x,y,k_1,k_2)|\nonumber\\
&\leq& \int_0^{+\infty} \left|\mathbb{E} \Big(D^2_{xx}(\bar{a}(x)-
a(x,Y_s^x(y)))\cdot
(k_1,k_2),D_{x}\bar{u}(t,x)\Big)_{\mathbb{R}^n }\right| ds\nonumber\\
&\leq&
C_T\left(1+\|x\|_{\mathbb{R}^n}+\|y\|_{\mathbb{R}^m}\right)\|k_1\|_{\mathbb{R}^n
}\|k_2\|_{\mathbb{R}^n
}\int_0^{+\infty}e^{-\frac{\beta}{2}s}ds\nonumber\\
&\leq&C_T\left(1+\|x\|_{\mathbb{R}^n}+\|y\|_{\mathbb{R}^m}\right)\|k_1\|_{\mathbb{R}^n
}\|k_2\|_{\mathbb{R}^n }.\label{J_1}
\end{eqnarray}
 By Lemma \ref{mix-derivative} and  \eqref{Averaging-Expectation} we infer that
\begin{eqnarray}
&&\!\!\!\!\!\!\!\!\!\!|J_2(t,x,y,k_1,k_2)|\nonumber\\
&\leq& \int_0^{{+\infty}}\left|\mathbb{E}\Big(D_x(\bar{ {a}}(x)-
a(x,Y_s^x(y)))\cdot k_1,D^2_{xx}\bar{u}(t,x)\cdot
k_2\Big)_{\mathbb{R}^n }\right|ds\nonumber\\
&\leq&
C_T\left(1+\|x\|_{\mathbb{R}^n}+\|y\|_{\mathbb{R}^m}\right)\|k_1\|_{\mathbb{R}^n
}\|k_2\|_{\mathbb{R}^n
}\int_0^{+\infty}e^{-\frac{\beta}{2}s}ds\nonumber\\
&\leq&C_T\left(1+\|x\|_{\mathbb{R}^n}+\|y\|_{\mathbb{R}^m}\right)\|k_1\|_{\mathbb{R}^n
}\|k_2\|_{\mathbb{R}^n }.\label{J_2}
\end{eqnarray}
With a similar argument we can also show that
\begin{eqnarray}
&&\!\!\!\!\!\!\!\!\!\!|J_3(t,x,y,k_1,k_2)|\nonumber\\
&\leq&C_T\left(1+\|x\|_{\mathbb{R}^n}+\|y\|_{\mathbb{R}^m}\right)\|k_1\|_{\mathbb{R}^n
}\|k_2\|_{\mathbb{R}^n }.\label{J_3}
\end{eqnarray}
By making use of Lemma \ref{bar-u-xxx} and
\eqref{Averaging-Expectation}, we get
\begin{eqnarray}
&&\!\!\!\!\!\!\!\!\!\!|J_4(t,x,y,k_1, k_2)|\nonumber\\
&\leq&C_T\|k_1\|_{\mathbb{R}^n}\cdot\|k_2\|_{\mathbb{R}^n}
\cdot(1+\|x\|_{\mathbb{R}^n}+\|y\|_{\mathbb{R}^m})\int_0^{+\infty} e^{-\frac{\beta}{2}s}ds\nonumber\\
&\leq&
C_T\|k_1\|_{\mathbb{R}^n}\cdot\|k_2\|_{\mathbb{R}^n}(1+\|x\|_{\mathbb{R}^n}+\|y\|_{\mathbb{R}^m}).
\label{J_4}
\end{eqnarray}
In view of the above estimates \eqref{J_1}, \eqref{J_2}, \eqref{J_3}
and \eqref{J_4}, we can conclude that there exists a constant $C_T$
such that
\begin{eqnarray*}
|D^2_{xx}u_1(t,x,y)\cdot (k_1, k_2)|\leq C_T
\|k_1\|_{\mathbb{R}^n}\cdot\|k_2\|_{\mathbb{R}^n}(1+\|x\|_{\mathbb{R}^n}+\|y\|_{\mathbb{R}^m}),\;t\in[0,T],
\end{eqnarray*}
which means that for fixed $y\in  \mathbb{R}^m$ and $t\in [0, T]$,
\begin{eqnarray*}
\|D^2_{xx}u_1(t,x,y) \|_{L( {\mathbb{R}^n},\mathbb{R})}\leq
C_T(1+\|x\|_{\mathbb{R}^n}+\|y\|_{\mathbb{R}^m}),
\end{eqnarray*}
where $\|\cdot\|_{L(\mathbb{R}^n,\mathbb{R})}$ denotes the usual
operator norm on Banach space consisting of  bounded and linear
operators from $\mathbb{R}^n$ to $\mathbb{R}$. As the diffusion
function $g$ is bounded, we get
\begin{eqnarray*}
Tr\Big(D^2_{xx}u_1(t,x,y)gg^T\Big)
&\leq& C_T\|D^2_{xx}u_1(t,x,y)\|_{L( {\mathbb{R}^n},\mathbb{R})} \nonumber\\
&\leq &C_T(1+\|x\|_{\mathbb{R}^n}+\|y\|_{\mathbb{R}^m}).
\label{step2}
\end{eqnarray*}\\
\textbf{Step 3:} Estimate of $\lambda_1[u_1(t,x+c(x
),y)-u_1(t,x,y)].$\\
By Lemma \ref{u-1-abso-lemma} and boundedness condition of $c(x)$,
we directly have
\begin{eqnarray*}
&&|\lambda_1[u_1(t,x+c(x ),y)-u_1(t,x,y)]|\nonumber\\
&&\lambda_1[|u_1(t,x+c(x ),y)|+|u_1(t,x,y)|]\nonumber\\
 &&\leq C_T(1
+\|x\|_{\mathbb{R}^n}+\|y\|_{\mathbb{R}^m}),\;t\in [0,
T].\label{step3}
\end{eqnarray*}
Finally, it is now easy to gather all previous estimates  for terms
in \eqref{L1u1} and conclude
\begin{eqnarray*}
\left|\mathcal {L}_1u_1(t,x,y)\right|\leq
C_T(1+\|x\|_{\mathbb{R}^n}+\|y\|_{\mathbb{R}^m}), \;t\in [0, T].
\end{eqnarray*}
\end{proof}
\end{lemma}

\begin{lemma}\label{r-bound-lemma}
Under the conditions of Lemma \ref{4.3}, for any $T>0$,
$x\in{\mathbb{R}^n}$ and $ y \in{\mathbb{R}^m}$, we have
\begin{eqnarray*}
|r^\epsilon(T,x,y)|\leq
C_T\epsilon(1+\|x\|_{\mathbb{R}^n}+\|y\|_{\mathbb{R}^m}).
\end{eqnarray*}
\begin{proof}
By a variation of constant formula, we  write the equation
\eqref{r-equation} in its integral form
\begin{eqnarray*}
r^\epsilon(T,x,y)&=&\mathbb{E}[r^\epsilon(0,X^\epsilon_T(x,y),Y^\epsilon_{T}(x, y))]\\
&+&\epsilon\left[\int_0^T\mathbb{E}(\mathcal{L}_1u_1-\frac{\partial
u_1}{\partial s})(s, X^\epsilon_{T-s}(x,y),Y^\epsilon_{ {T-s} }(x,
y)) ds\right].
\end{eqnarray*}
Since $u^\epsilon$ and $\bar{u}$ satisfy the same initial condition,
we have
\begin{eqnarray*}
|r^\epsilon(0, x,y)|&=&|u^\epsilon(0,x,y)-\bar{u}(0,x)-\epsilon
u_1(0,x,y)|\\
&=&\epsilon |u_1(0,x,y)|,
\end{eqnarray*}
so that, thanks to \eqref{u-1-abso}, \eqref{X-bound} and
\eqref{Y-bound} we have
\begin{eqnarray}
\mathbb{E}[r^\epsilon(0,X^\epsilon_T(x,y),Y^\epsilon_{T }(x,y)]\leq
C\epsilon(1+\|x\|_{\mathbb{R}^n}+ \|y\|_{\mathbb{R}^m}).
\label{r-bound}
\end{eqnarray}
Using Lemma \ref{4.3} and Lemma \ref{L2u1} yields
\begin{eqnarray*}
&&\mathbb{E}[(\mathcal{L}_1u_1-\frac{\partial u_1}{\partial s})(s,
X^\epsilon_{T-s}(x,y),Y^\epsilon_{T-s} (x,y))]\\
&&\leq C \mathbb{E}
(1+\|X^\epsilon_{T-s}(x,y)\|+\|Y^\epsilon_{T-s}(x,y)\|),
\end{eqnarray*}
and, according to \eqref{X-bound} and \eqref{Y-bound}, this implies
that
\begin{eqnarray*}
&&\mathbb{E}\left[\int_0^T(\mathcal{L}_1u_1-\frac{\partial
u_1}{\partial s})(s, X^\epsilon_{T-s}(x,y), Y^\epsilon_{T-s}(x,y))
ds\right]\nonumber\\
&&\leq C_T (1+\|x\|_{\mathbb{R}^n}+ \|y\|_{\mathbb{R}^m}).
\end{eqnarray*}
The last inequality together with \eqref{r-bound}  yields
\begin{eqnarray*}
|r^\epsilon(T,x,y)|\leq  \epsilon C_T (1+\|x\|_{\mathbb{R}^n}+
\|y\|_{\mathbb{R}^m}).
\end{eqnarray*}
\end{proof}
\end{lemma}

\section{Appendix}\label{appendix}
In this appendix we collect  some technical results to which we
appeal in the proofs of the main results in Section 4 .
\begin{lemma}\label{bar-u-x}
For any $T>0$, there exists a constant $C_T>0$ such that for any
$x,k\in\mathbb{R}^n$ and $t\in [0, T]$, we have
$$|D_{x}\bar{u}(t,x)\cdot k| \leq C_{T}\|k\|_{\mathbb{R}^n}.$$
\begin{proof}
Observe that for any   $k\in {\mathbb{R}^n}$,
\begin{eqnarray*}
D_x\bar{u}(t,x)\cdot
k&=&\mathbb{E}\left[D\phi(\bar{X}_t(x))\cdot\eta ^{k,x}_t\right]\\
&=&\mathbb{E}\left(\phi'(\bar{X}_t(x)),\eta^{k,x}_t\right)_{\mathbb{R}^n},
\end{eqnarray*}
where $\eta^{k,x}_t$ denotes the first mean-square derivative of
$\bar{X}_t(x)$ with respect to $x\in \mathbb{R}^n$ along the
direction $k\in \mathbb{R}^n$, then we have
\begin{eqnarray*}
\begin{cases}
d\eta^{k,x}_t=D_x\bar{a}(\bar{X}_t(x))\cdot\eta^{k,x}_tdt+D_xb(\bar{X}_t(x))\cdot\eta^{k,x}_td {B}_t\\
\qquad\quad+D_xc(\bar{X}_{t-}(x))\cdot\eta^{k,x}_{t-}d {P}_t,\\
\eta ^{k,x}_0=k.
\end{cases}
\end{eqnarray*}
This means that  $\eta ^{k, x}_t$  is the  solution of the integral
equation
 \begin{eqnarray*}
\eta^{k, x}_t&=&k+\int_0^tD_x\bar{a}(\bar{X}_s(x))\cdot\eta^{k, x}_s
ds+\int_0^tD_xb(\bar{X}_s(x))\cdot\eta^{k, x}_sdB_s\\
&&+\int_0^tD_xc(\bar{X}_{s-}(x))\cdot \eta_{s-}^{k, x}dP_s
\end{eqnarray*}
and then thanks to  assumption (\textbf{A1}), we get
\begin{eqnarray*}
\mathbb{E}\|\eta ^{k,x}_t\|^2_{\mathbb{R}^n} \leq
C_T\|k\|^2_{\mathbb{R}^n} +C_T\int_0^t\mathbb{E}\|\eta
^{k,x}_s\|^2_{\mathbb{R}^n} ds.
\end{eqnarray*}
Then by  Gronwall lemma it follows that
\begin{eqnarray}
\mathbb{E}\|\eta ^{k,x}_t\|^2_{\mathbb{R}^n} \leq C_T\|k
\|^2_{\mathbb{R}^n}, \;t\in [0, T],  \label{eta-bound}
\end{eqnarray}
so that
\begin{eqnarray*}
|D_x\bar{u}(t,x)\cdot k|\leq  C_T \|k\|_{\mathbb{R}^n}.
\end{eqnarray*}
\end{proof}
\end{lemma}
Next, we introduce an analogous result for the second derivative of
$\bar{u}(t,x)$.
\begin{lemma}\label{bar-u-xx}
For any $T>0$, there exists a constant $C_T>0$ such that for any
$x,k_1,k_2\in\mathbb{R}^n$ and $t\in [0, T]$, we have
$$|D^2_{xx}\bar{u}(t,x)\cdot(k_1,k_2)|\leq C_{T}\|k_1\|_{\mathbb{R}^n}\cdot\|k_2\|_{\mathbb{R}^n}.$$
\begin{proof}
For any $k_1, k_2 \in \mathbb{R}^n$, we have
\begin{eqnarray}
D^2_{xx}\bar{u}(t,x)\cdot(k_1,k_2)&=&\mathbb{E}\big[\phi''(\bar{X}_t(x))\cdot(\eta^{k_1,x}_t,\eta^{k_2,x}_t)\nonumber\\
&+&\phi'(\bar{X}_t(x))\cdot \xi^{k_1,k_2,x}_t\big],
\label{Dxx-bar-u}
\end{eqnarray}
where $\xi^{k_1,k_2,x}_t$ is  the solution of     the second
variation equation corresponding to the averaged equation, which may
be rewritten in the following form:
\begin{eqnarray*}
\xi^{k_1,k_2,x}_t&=&\int_0^t\big[D_x\bar{a} (\bar{X}_s(x))\cdot\xi^{k_1,k_2,x}_s+D_{xx}^2\bar{a} (\bar{X}_s(x))\cdot(\eta^{k_1,x}_s,\eta^{k_2,x}_s)\big]ds\\
&+& \int_0^t\big[D_{xx}^2b (\bar{X}_s(x))\cdot(\eta^{k_1,x}_s,\eta^{k_2,x}_s)+D_xb(\bar{X}_s(x))\cdot\xi^{k_1,k_2,x}_s\big]dB_s\\
&+&\int_0^t\big[D_{xx}^2c(\bar{X}_{s-}(x))\cdot(\eta^{k_1,x}_{s-},\eta^{k_2,x}_{s-})+D_xc(\bar{X}_{s-}(x))\cdot\xi^{k_1,k_2,x}_{s-}\big]dP_s.
\end{eqnarray*}
Thus, by assumption (\textbf{A1}) and \eqref{eta-bound}  we have
\begin{eqnarray*}
\mathbb{E}\|\xi^{k_1,k_2,x}_t\|^2_{\mathbb{R}^n}&\leq&
C_T\int_0^t\big(\{\mathbb{E}\|\eta^{k_1,x}_s\|^2_{\mathbb{R}^n}\}^{\frac{1}{2}}
\{\mathbb{E}\|\eta^{k_2,x}_s\|^2_{\mathbb{R}^n}\}^{\frac{1}{2}}+\mathbb{E}\|\xi^{k_1,k_2,x}_s\|^2_{\mathbb{R}^n}\big)ds\\
&\leq& C_T\|k_1\|_{\mathbb{R}^n}\cdot\|k_2\|_{\mathbb{R}^n}
+C_T\int_0^t\mathbb{E}\|\xi^{k_1,k_2,x}_s\|^2_{\mathbb{R}^n}ds.
\end{eqnarray*}
By applying the Gronwall lemma we have
\begin{equation*}
\mathbb{E}\|\zeta^{k_1,k_2,x}_t\|^2_{\mathbb{R}^n}\leq
C_T\|k_1\|_{\mathbb{R}^n}\cdot \|k_2\|_{\mathbb{R}^n}.
\end{equation*}
Returning to \eqref{Dxx-bar-u}, we can get
\begin{eqnarray*}
|D^2_{xx}\bar{u}(t,x)\cdot(k_1,k_2)|\leq C_T
\|h_1\|_{\mathbb{R}^n}\cdot\|k_2\|_{\mathbb{R}^n}.
\end{eqnarray*}
\end{proof}
\end{lemma}
By using the analogous arguments used before, we can prove the
following estimate for the  third order derivative of $\bar{u}(t,x)$
with respect to $x$.
\begin{lemma}\label{bar-u-xxx}
For any $T>0$, there exists a constant $C_T>0$ such that for any
$x,k_1,k_2,k_3\in\mathbb{R}^n$ and $t\in [0, T]$, we have
$$|D^3_{xxx}\bar{u}(t,x)\cdot(k_1,k_2,k_3)|\leq C_{T }\|k_1\|_{\mathbb{R}^n}\cdot\|k_2\|_{\mathbb{R}^n}\cdot\|k_3\|_{\mathbb{R}^n}.$$
\end{lemma}

The following lemma states   boundedness for the first derivative of
$\bar{a}(x)-\mathbb{E}a(x, Y^x_t(y))$ with respect to $x$.
\begin{lemma}\label{mix-derivative}
There exists a constant $C>0$ such that for any $x\in \mathbb{R}^n,
y\in \mathbb{R}^m$, $k\in \mathbb{R}^n$ and $t>0$ it holds
\begin{eqnarray*}
\|D_x (\bar{a}(x)-\mathbb{E}a(x, Y^x_t(y)))\cdot
k\|_{\mathbb{R}^n}\leq
Ce^{-\frac{\beta}{2}t}\|k\|_{\mathbb{R}^n}\left(1+\|x\|_{\mathbb{R}^n}+\|y\|_{\mathbb{R}^m}\right).
\end{eqnarray*}
\begin{proof}
 The proof is a modification of the proof of \cite[Proposition C.2]{Brehier}. For
any $t_0>0$, we set
\begin{eqnarray*}
\tilde{a}_{t_0}(x,y,t)=\hat{a}(x,y,t)-\hat{a}(x,y,t+t_0),
\end{eqnarray*}
where
\begin{eqnarray*}
\hat{a}(x,y,t):=\mathbb{E}a(x, Y^{x}_t(y)).
\end{eqnarray*}
Then we have
\begin{eqnarray*}
\lim\limits_{t_0\rightarrow
+\infty}\tilde{a}_{t_0}(x,y,t)=\mathbb{E}a(x, Y^x_t(y))-\bar{a}(x).
\end{eqnarray*}
By Markov property, we have
\begin{eqnarray*}
\tilde{a}_{t_0}(x,y,t)&=&\hat{a}(x,y,t)-\mathbb{E}{a}(x,Y_{t+t_0}^{x}(y))\\
&=&\hat{a}(x,y,t)-\mathbb{E}\hat{a}(x, Y_{t_0}^{x}(y),t)
\end{eqnarray*}
Due to assumption (\textbf{A1}),  for any $k\in \mathbb{R}^n$ we
have
\begin{eqnarray}
D_x\tilde{a}_{t_0}(x,y,t)\cdot k&=&D_x\hat{a}(x,y,t)\cdot
k-\mathbb{E}D_x\left(\hat{a}(x, Y_{t_0}^{x }(y),t)\right)\cdot k\nonumber\\
&=&\hat{a}_x'(x,y,t)\cdot k-\mathbb{E}\hat{a}_x'(x,
Y_{t_0}^{x}(y),t)\cdot
k\nonumber\\
&&-\mathbb{E}\hat{a}_y'(x,
Y_{t_0}^{x}(y),t)\cdot\left(D_xY_{t_0}^{x}(y)\cdot
k\right),\label{7-1-1}
\end{eqnarray}
where  the symbols $\hat{a}_x'$ and $\hat{a}_y'$   denote the
directional derivatives with respect to $x$ and $y$, respectively.
Note that the first derivative $\zeta_t^{x,y,
k}=D_xY_{t}^{x}(y)\cdot k$, at the point $x$ and along the direction
$k\in \mathbb{R}^n$, is the solution of equation
\begin{eqnarray*}
d\zeta_t^{x, y,k}&=&\left( f_x'(x, Y_t^{x}(y))\cdot
k+f_y'(x, Y_t^{x}(y))\cdot\zeta_t^{x, y,k}\right)dt\\
&&+\left(g_x'(x, Y_t^{x}(y))\cdot k+g_y'(x,
Y_t^{x}(y))\cdot\zeta_t^{x, y,k}\right)dW_t\\
&&+\left(h_x'(x, Y_{t-}^{x}(y))\cdot k+h'_y(x,
Y_{t-}^{x}(y))\cdot\zeta_{t-}^{x, y,k}\right)dN_t
\end{eqnarray*}
with initial data $\zeta_0^{x,y, k}=0$. Hence,   by assumption
(\textbf{A1}), it is straightforward to check
\begin{eqnarray}
\mathbb{E}\|\zeta_t^{x,y, k}\|_{\mathbb{R}^m}\leq
C\|k\|_{\mathbb{R}^n}\label{7-2}
\end{eqnarray}
for any $t\geq 0$. Note that for any $y_1, y_2\in {\mathbb{R}^m}$,
we have
\begin{eqnarray}
\|\hat{a}(x,y_1,
t)-\hat{a}(x,y_2,t)\|_{\mathbb{R}^n}&=&\|\mathbb{E}a(x,
Y_t^{x}(y_1))-\mathbb{E}a(x,
Y_t^{x}(y_2))\|_{\mathbb{R}^n}\nonumber\\
&\leq&C\mathbb{E}\|Y_t^{x}(y_1)-Y_t^{x}(y_2)\|_{\mathbb{R}^m}\nonumber\\
&\leq& Ce^{-\frac{\beta}{2}t}\|y_1-y_2\|_{\mathbb{R}^m},\nonumber
\end{eqnarray}
where   \eqref{difference-initial} was used to obtain the last
inequality. This means that
\begin{eqnarray}
\|\hat{a}_y'(x, y,t)\cdot l\|_{\mathbb{R}^m}\leq
Ce^{-\frac{\beta}{2}t}\|l\|_{\mathbb{R}^m},\; l\in
{\mathbb{R}^m}.\label{7-3}
\end{eqnarray}
From  \eqref{7-2} and \eqref{7-3}, we obtain
\begin{eqnarray}
&&\|\mathbb{E}[\hat{a}_y'(x,
Y_{t_0}^{x}(y),t)\cdot\left(D_xY_{t_0}^{x}(y)\cdot
k\right)]\|_{\mathbb{R}^m}\nonumber\\
&&=\|\mathbb{E}[\hat{a}_y'(x,
Y_{t_0}^{x}(y),t)\cdot (\zeta_{t_0}^{x,y, k} )]\|_{\mathbb{R}^m}\nonumber\\
&&\leq C e^{-\frac{\beta}{2}t}\|k\|_{\mathbb{R}^n}.\label{7-4}
\end{eqnarray}
Then, by easy calculations, we have
\begin{eqnarray}
&&\hat{a}_x'(x,y_1,t)\cdot k-\hat{a}_x'(x,y_2,t)\cdot k\nonumber\\
&&\quad=\mathbb{E}\left(a_x'(x, Y_t^{x}(y_1))\right)\cdot
k-\mathbb{E}\left(a_x'(x, Y_t^{x}(y_2))\right)\cdot k\nonumber\\
&&\quad\quad+\mathbb{E}\left(a_y'(x, Y_t^{x}(y_1))\cdot
\zeta_t^{x,y_1,
k}-a_y'(x, Y_t^{x}(y_2))\cdot \zeta_t^{x,y_2, k}\right)\nonumber\\
&&\quad= \mathbb{E}\left(a_x'(x, Y_t^{x}(y_1))\right)\cdot
k-\mathbb{E}\left(a_x'(x, Y_t^{x}(y_2))\right)\cdot k\nonumber\\
&&\quad\quad+\mathbb{E}\left([a_y'(x, Y_t^{x}(y_1))-a_y'(x,
Y_t^{x}(y_2))]\cdot\zeta_t^{x,y_1, k} \right)\nonumber\\
&&\quad\quad+\mathbb{E}\left(a_y'(x,
Y_t^{x}(y_2))\cdot(\zeta_t^{x,y_1, k}-\zeta_t^{x,y_2,
k})\right)\nonumber \\
&&\quad:= \sum\limits_{i=1}^3\mathcal {N}_i(t,x,y_1,y_2, k).
\label{7-5}
\end{eqnarray}
Now, we estimate the three terms in the right hand side of above
equality. Concerning $\mathcal{N}_1(t,x,y_1,y_2, k)$ we have
\begin{eqnarray}
&&\|\mathcal{N}_1(t,x,y_1,y_2, k)\|_{\mathbb{R}^n}\nonumber\\
&& \leq\mathbb{E}\|\left(a_x'(x, Y_t^{x}(y_1))\right)\cdot
k-\left(a_x'(x, Y_t^{x}(y_2))\right)\cdot k\|_{\mathbb{R}^n}\nonumber\\
&& \leq
C\mathbb{E}\|Y_t^{x}(y_1)-Y_t^{x}(y_2)\|_{\mathbb{R}^m}\cdot\|k\|_{\mathbb{R}^n}\nonumber\\
&& \leq
Ce^{-\frac{\beta}{2}t}\|y_1-y_2\|_{\mathbb{R}^m}\cdot\|k\|_{\mathbb{R}^n}.\label{7-6}
\end{eqnarray}
Next, by  assumption (\textbf{A1})  we get
\begin{eqnarray}
&&\|\mathcal{N}_2(t,x,y_1,y_2, k)\|_{\mathbb{R}^n}\nonumber\\
&&\leq\mathbb{E}\|[a_y'(x, Y_t^{x}(y_1))-a_y'(x,
Y_t^{x}(y_2))]\cdot\zeta_t^{x,y_1, k}\|_{\mathbb{R}^n}\nonumber\\
&&\leq C\{\mathbb{E}\|\zeta_t^{x,y_1,
k}\|^2_{\mathbb{R}^m}\}^{\frac{1}{2}}\cdot\{\mathbb{E}\|Y_t^{x}(y_1)-Y_t^{x}(y_2)\|^2_{\mathbb{R}^m}\}^{\frac{1}{2}}\nonumber\\
&&\leq C
e^{-\frac{\beta}{2}t}\|k\|_{\mathbb{R}^n}\cdot\|y_1-y_2\|_{\mathbb{R}^m}.\label{7-7}
\end{eqnarray}
For the third term, by making use of assumption (\textbf{A1}) again,
we can infer that
\begin{eqnarray}
&&\| \mathcal{N}_3(t,x,y_1,y_2, k)\|_{\mathbb{R}^n}\nonumber\\
&& \leq\mathbb{E}\|a_y'(x, Y_t^{x}(y_2))\cdot(\zeta_t^{x,y_1, k}-\zeta_t^{x,y_2, k})\|_{\mathbb{R}^n}\nonumber\\
&& \leq C\mathbb{E}\|\zeta_t^{x,y_1, k}-\zeta_t^{x,y_2, k}\|_{\mathbb{R}^m}\nonumber\\
&& \leq C
e^{-\frac{\beta}{2}t}\|y_1-y_2\|_{\mathbb{R}^m}\cdot\|k\|_{\mathbb{R}^n}.\label{7-8}
\end{eqnarray}
Now, returning to \eqref{7-5} and  taking into account of
\eqref{7-6}, \eqref{7-7} and \eqref{7-8}, we get
\begin{eqnarray*}
&&\|\hat{a}_x'(x,y_1,t)\cdot k-\hat{a}_x'(x,y_2,t)\cdot k\|\nonumber\\
&&\leq C
e^{-\frac{\beta}{2}t}\|y_1-y_2\|_{\mathbb{R}^m}\cdot\|k\|_{\mathbb{R}^n},
\end{eqnarray*}
which leads to
\begin{eqnarray}
&&\|\hat{a}_x'(x,y,t)\cdot h- \mathbb{E}\hat{a}_x'(x,Y_{t_0}^{x}(y),t)\cdot k\|_{\mathbb{R}^n}\nonumber\\
&&\leq C
e^{-\frac{\beta}{2}t}(1+\|y\|_{\mathbb{R}^m}+\|Y_{t_0}^{x}(y)\|_{\mathbb{R}^m})\cdot\|k\|_{\mathbb{R}^n}\nonumber\\
&&\leq
e^{-\frac{\beta}{2}t}(1+\|x\|_{\mathbb{R}^n}+\|y\|_{\mathbb{R}^m})\cdot\|k\|_{\mathbb{R}^n},\label{7-9}
\end{eqnarray}
where we used the inequality \eqref{frozen-bound}. Returning to
\eqref{7-1-1}, by \eqref{7-4} and \eqref{7-9} we conclude that
\begin{eqnarray*}
\|D_x\tilde{a}_{t_0}(x,y,t)\cdot k\|_{\mathbb{R}^n}\leq
Ce^{-\frac{\beta}{2}t}(1+\|x\|_{\mathbb{R}^n}+\|y\|_{\mathbb{R}^m})\|k\|_{\mathbb{R}^n}.
\end{eqnarray*}
Taking the limit as $t_0\rightarrow +\infty$  we obtain
\begin{eqnarray*}
\|D_x (\bar{a}(x)-\mathbb{E}a(x, Y^x_t(y)))\|_{\mathbb{R}^n}\leq
Ce^{-\frac{\beta}{2}{t}}\|k\|_{\mathbb{R}^n}\left(1+\|x\|_{\mathbb{R}^n}+\|y\|_{\mathbb{R}^m}\right).
\end{eqnarray*}

\end{proof}
\end{lemma}
 Proceeding with similar arguments above we obtain the following higher order
 differentiability.
\begin{lemma}\label{mix-derivative-2}
There exists a constant $C>0$ such that for any $x, k_1, k_2\in
\mathbb{R}^n$, $y\in \mathbb{R}^m$ and $t>0$ it holds
\begin{eqnarray*}
&&\|D^2_{xx}(\bar{a}(x)-\mathbb{E}a(x,
Y^x_t(y)))(k_1,k_2)\|_{\mathbb{R}^n}\\
&&\leq
Ce^{-\frac{\beta}{2}{t}}\|k_1\|_{\mathbb{R}^n}\|k_2\|_{\mathbb{R}^n}\left(1+\|x\|_{\mathbb{R}^n}+\|y\|_{\mathbb{R}^m}\right).
\end{eqnarray*}
\end{lemma}

Finally, we introduce the following auxiliary result.
\begin{lemma}\label{deri-xt-bar-u}
There exists a constant $C>0$ such that for any $x, k\in
\mathbb{R}^n$, $y\in \mathbb{R}^m$ and $t>0$ it holds
$$\|\frac{\partial}{\partial t}D_x\bar{u}(t,x)\cdot k\|_{\mathbb{R}^n}\leq C\|k\|_{\mathbb{R}^n}.$$
\begin{proof}
For simplicity of presentation, we will prove it for the
1-dimensional case. The multi-dimensional situation can be treated
similarly, only notations are somewhat involved.  In this case we
only need to show
\begin{eqnarray}
|\frac{\partial}{\partial t}\frac{\partial}{\partial
x}\bar{u}(t,x)|\leq C.\label{1-dimen}
\end{eqnarray}
Actually, for any $\phi\in C_b^3(\mathbb{R},\mathbb{R})$ we have
\begin{eqnarray*}
\frac{\partial}{\partial x}\bar{u}(t,x)=\frac{\partial}{\partial
x}\mathbb{E}\phi(\bar{X}_t(x))=\mathbb{E}\left(\phi'(\bar{X}_t(x))\cdot\frac{\partial}{\partial
x}\bar{X}_t(x)\right).
\end{eqnarray*}
If we define
\begin{eqnarray*}
\varsigma^x_t:=\frac{\partial}{\partial x}\bar{X}_t(x),
\end{eqnarray*}
we have
\begin{eqnarray*}
\varsigma^x_t&=&1+\int_0^t\bar{a}'(\bar{X}_s(x))\cdot\varsigma^x_s
ds+\int_0^tb'(\bar{X}_s(x))\cdot\varsigma^x_sdB_s\\
&&+\int_0^tc'(\bar{X}_{s-}(x))\cdot \varsigma_{s-}^xdP_s.
\end{eqnarray*}
The boundedness of $\bar{a}', b'$ and $c'$   guarantees
\begin{eqnarray}
\mathbb{E}|\varsigma^x_t|^2\leq C_T, \;t\in [0,
T].\label{varsigma-bound}
\end{eqnarray}
 By    using It\^{o} formula we have
\begin{eqnarray*}
&&\mathbb{E}[\phi'(\bar{X}_t(x))\cdot\varsigma_t^x]\\
&&=\phi'(x)+\mathbb{E}\int_0^t[\phi'(\bar{X}_s(x))\bar{a}'(\bar{X}_s(x))\varsigma_s^x
+\varsigma_s^x\phi''(\bar{X}_s(x))\bar{a}(\bar{X}_s(x))]ds\\
&&+\mathbb{E}\int_0^tb'(\bar{X}_s(x)))\varsigma_s^x\phi''(\bar{X}_s(x))b(\bar{X}_s(x))ds\\
&&+\frac{1}{2}\mathbb{E}\int_0^t\varsigma_s^x\phi'''(\bar{X}_s(x))b^2(\bar{X}_s(x))ds\\
&&+\lambda_1\mathbb{E}\int_0^t\phi'(\bar{X}_s(x))c'(\bar{X}_{s-}(x))\varsigma_s^xds\\
&&+\lambda_1\mathbb{E}\int_0^t\varsigma_s^x[\phi'(\bar{X}_{s-}(x)+c(\bar{X}_{s-}(x)))-\phi'(\bar{X}_{s-}(x))]ds\\
&&+\lambda_1\mathbb{E}\int_0^tc'(\bar{X}_{s-}(x))\varsigma_s^x[\phi'(\bar{X}_{s-}(x)+c(\bar{X}_{s-}(x)))-\phi'(\bar{X}_{s-}(x))]ds.
\end{eqnarray*}
 Since $\phi$ belongs to $C_b^3(\mathbb{R}, \mathbb{R})$, from
 the assumption (\textbf{A1}) it follows that for any $t\in [0, T],$
\begin{eqnarray*}
\left|\frac{\partial}{\partial t}[\frac{\partial}{\partial
x}\bar{u}(t,x)]\right|&=&\left|\frac{\partial}{\partial t}\mathbb{E}[\phi'(\bar{X}_t(x))\cdot\varsigma_t^x]\right|\\
&\leq& C|\mathbb{E}\varsigma_t^x|,
\end{eqnarray*}
then, by taking \eqref{varsigma-bound} into account, one would
easily arrive at \eqref{1-dimen}.
\end{proof}
\end{lemma}

\section*{Acknowledgments}
{We would like to thank Professor} Jinqiao Duan for helpful
discussions and comments. Hongbo Fu is supported by Natural Science
Foundation of Hubei Province (No. 2018CFB688), NSF of China (No.
11301403) and  Chinese Scholarship Council (No. [2015]5104). Bengong
Zhang is supported by NSF of China (No. 11401448).  Li Wan is
supported by NSF  of China (No. 61573011). Jicheng Liu is supported
by NSF  of China (No. 11271013).

\section*{Competing interests}
The authors declare  that they have no competing interests.

\section*{Author¡¯s contributions}
The authors declare that the work was realized in collaboration with
the same responsibility. All authors read and approved the final
manuscript.

\label{}











\end{document}